\documentclass{amsart}
\usepackage{amsmath}
\usepackage{graphicx} 
\newtheorem{theorem}{Theorem}
\newtheorem{conjecture}[theorem]{Conjecture}

\newtheorem{remarks}[theorem]{Remarks}
\newtheorem{example}[theorem]{Example}

\usepackage{hyperref,xcolor}
\usepackage[T1]{fontenc}
\usepackage{pdfsync}
\usepackage{yfonts}      
\usepackage{tcolorbox}   
\usepackage{lipsum} 
\DeclareMathSymbol{\lsim}{\mathord}{symbols}{"18}
\title{Using Walnut to solve problems from the OEIS}
\author[Bosma\and Bruin\and Fokkink\and Grube\and Reuijl\and Tromp]{Wieb Bosma\and Ren\'e Bruin\and Robbert Fokkink\and Jonathan Grube\and Anniek Reuijl\and Thian Tromp}
\email{bosma@math.ru.nl \and r.j.fokkink@tudelft.nl} 

\begin{document}
\begin{abstract}
We use the automatic theorem prover \texttt{Walnut} to resolve various open problems from the OEIS and beyond. Specifically, we clarify the structure of sequence A260311, which concerns runs of sums of upper Wythoff numbers. We extend a result of Hajdu, Tijdeman, and Varga on polynomials with nonzero coefficients modulo a prime. Additionally, we settle open problems related to the anti-recurrence sequences A265389 and A299409, as well as the subsumfree sequences A026471 and A026475. Our findings also give rise to new open problems.
\end{abstract}

\subjclass[2000]{11B37, 11R09, 68R15}
\maketitle

In the fall of 2024 two of us (Bosma and Fokkink) taught a course on
automatic sequences
in the Dutch national mathematics master's programs, based on Jeffrey Shallit's 
monograph~\cite{shallit2022}. The aim
of the course was to allow students to solve as many problems as 
possible from the Online Encyclopedia of Integer Sequences using
the automatic theorem prover \texttt{Walnut}~\cite{mousavi2016}. This paper contains
our results. It consists of four distinct sections: 
an analysis of sums of upper Wythoff numbers, an investigation of polynomials with rational roots, 
and two separate studies on specific sum-free sequences.
All these studies are carried out with~\texttt{Walnut}.

\texttt{Walnut} is free software that is available from
\href{https://cs.uwaterloo.ca/~shallit/walnut.html}{\small{cs.uwaterloo.ca/\~{}shallit/walnut.html}}.
It can verify first-order logic statements on sets or sequences of natural numbers involving addition but
not multiplication. Statements can be expressed in various numeration systems, such as binary,
decimal, or the exotic Zeckendorf numeration system, in which numbers are represented by sums
of Fibonacci numbers.
The syntax of \texttt{Walnut} is highly transparent and can be easily understood by
anyone familiar with first-order logic. To illustrate this, here is
an example. 

The Fibonacci word $\mathbf F=abaababaabaab\cdots$
is the infinite symbolic sequence that is generated by the substitution $a\mapsto ab$ and $b\mapsto a$.
A well-known property of $\mathbf F$ is that each prefix of length $F_n-2$ is a palindrome~\cite{deluca1981},
where $F_n$ is the $n$-th Fibonacci number. 
Such a prefix runs up to index $F_n-3$ if we start the count at zero.
In first order logic, the palindromic property can be expressed as
\[
\forall F_n\geq 3\ \forall j\leq F_n-3\colon\ \mathbf F[j]=\mathbf F[F_n-3-j].
\]
The Fibonacci numbers $0,1,1,2\ldots$ are denoted $F_0,F_1,F_2,F_3,\ldots$ and
therefore $F_4$ is the first Fibonacci number such that $F_n\geq 3$. A \texttt{Walnut} verification of the statement is
\begin{center}
    \texttt{eval prefixtest "?msd\_fib Af Aj \\ (\$isfib(f) \& f>=3 \& j<=f-3) => F[j]=F[f-3-j]":}
\end{center}

We explain the command.
It determines the validity of the hyphenated statement, which
is named \texttt{prefixtest}. The command \texttt{?msd\_fib} indicates that it is a
statement in the Zeckendorf system.
The abbreviation \texttt{msd} stands for most significant digit first.
The capital A represents $\forall$, and similarly, E represents $\exists$.
The command \texttt{\$isfib(f)} checks whether \texttt{f} is a Fibonacci
number. 
The Fibonacci word $\mathbf F$ is implemented in \texttt{Walnut} as~\texttt{F}.
The command should now be clear. It is evaluated as \texttt{TRUE}.

\section{Gaps between unsums of upper Wythoff numbers}

The upper Wythoff numbers are the indices for the letter $b$ in the
word $\mathbf F$, if we start the count at one instead of zero. The
first few upper Wythoff numbers are
\[
2,\ 5,\ 7,\ 10,\ 13,\ 15,\ 18,\ 20,\ 23,\ 26,\ 28,\ 31,\ \ldots
\]
Sequence \href{https://www.oeis.org/a260317}{A260317} in the OEIS 
contains the
numbers that are not sums of upper Wythoff numbers. 
We say that they are \emph{unsums}.
Sequence \href{https://www.oeis.org/a260311}{A260311} in the OEIS contains the
gaps between consecutive unsums.
The initial entries of this sequence are
{\small
\[
\textcolor{gray}{1111122121}\textcolor{red}{2}3\textcolor{red}{2}3\textcolor{red}{2}33\textcolor{red}{2}33\textcolor{red}{2}353\textcolor{red}{2}353\textcolor{red}{2}35353\textcolor{red}{2}35353\textcolor{red}{2}35355353\textcolor{red}{2}35355353\textcolor{red}{2}3535535355353\textcolor{red}{2}
\]
}

\noindent Jeffrey Shallit~\cite{shallit2023}
implemented \texttt{a260317} in \texttt{Walnut} to prove that
all entries in  \href{https://www.oeis.org/a260311}{A260311} are limited to 1, 2, 3, and 5.
This resolved a problem posed by Clark Kimberling.
Building on this work, we extend the analysis of \href{https://www.oeis.org/a260311}{A260311}
and fully clarify its structure.

The prefix $1111122121$ is irregular, but after that a pattern emerges. 
We first describe what we see and then we check that with \texttt{Walnut}.
The words between
the entries $2$ are palindromes of increasing lengths:
\[
\textcolor{red}{2}W_1\textcolor{red}{2}W_1\textcolor{red}{2}W_2\textcolor{red}{2}W_2\textcolor{red}{2}W_3\textcolor{red}{2}W_3\textcolor{red}{2}\cdots
\]
Each palindrome $W_i$ occurs twice and has Fibonacci length $|W_i|=F_{i+1}$.
The sequence of palindromes $W_1, W_2, W_3,\ldots$ is 
\[3,\ 33,\ 353,\ 35353,\ 35355353,\ 3535535355353,\ 353553535535535355353,\ \ldots \]
From the third palindrome $W_3$ on, each next palindrome $W_{n+1}$ agrees with
$W_n$ up to the its $k$-th letter, with $k=\lceil|W_{n+1}|/2\rceil$. 
In particular, each $W_n$ determines $W_{n+1}$ since a palindrome is determined
by its initial half: all palindromes follow from $W_3$, since we know their lengths.
The digit sums of the palindromes are
\[3,\ 6,\ 11,\ 19,\ 32,\ 53,\ 87,\ \ldots\]
and if we add $2$ to these sums, which is the preceding digit, 
then we get the Fibonacci numbers $F_5,F_6,F_7,F_8 \ldots$. The limit of
the sequence of palindromes is the infinite word
\[
3535535355355353553535535535355355\cdots.
\]
If we replace $a$ by $5$ and $b$ by $3$ in the Fibonacci word \textbf{F}, then this
limit is equal to~$3\mathbf{F}$.
All these properties can be checked with \texttt{Walnut}.

We first make sure that gap $1$ only occurs in the
the prefix $1111122121$. 
The digit sum of the prefix is equal to $13$ and the first
entry of \href{https://www.oeis.org/a260311}{A260311} is $1$.
That is why that last gap of $1$ in the prefix is the gap between
$13$ and $14$ in \href{https://www.oeis.org/a260317}{A260317}.
We verify that the final gap $1$ indeed occurs between the entries $13$ and $14$ in this sequence:
\begin{center}
    \texttt{eval test "?msd\_fib An (\$a260317(n) \& \$a260317(n+1)) => n<14":}
\end{center}
\texttt{Walnut} returns \texttt{TRUE}.
The digits $\textcolor{red}{2}$ are markers and we single them out. Since there
are no gaps $1$ from number $14$ on, all gaps are $\geq 2$ for $n>13$:
\begin{center}
    \texttt{def marker2 "?msd\_fib \$a260317(n) \& \$a260317(n+2) \& n>13":}
\end{center}

\noindent We can now check that the digit sums of the palindromes
plus two are the Fibonacci numbers starting from $5$. Since the digit sum
is a sum of gaps, it is equal to the difference between
consecutive markers. We first collect this difference in
\texttt{gapmark2} and then check that these are the Fibonacci numbers
starting from $F_5$.
\begin{center}
    \texttt{def gapmark2 "?msd\_fib En (t>0) \& \$marker2(n) \& \$marker2(n+t) \&}
    
    \texttt{(As (s<t \& s>0) => \~{}\$marker2(n+s))":}

\texttt{eval test "?msd\_fib An \$gapmark2(n) <=> (n>4 \& \$isfib(n))":}
\end{center}

\noindent\texttt{Walnut} returns \texttt{TRUE}. We already saw that the digit sums
are $F_5, F_5, F_6, F_6, \ldots$. 
Consecutive triples in this sequence are
either $F_k,F_k,F_{k+1}$ or $F_k,F_{k+1},F_{k+1}$ for $k>4$. 
To check this with \texttt{Walnut} we use that 
these are
exactly the triples $(t,u,v)$ of non-decreasing Fibonacci numbers such that
$t+v$ is a Fibonacci number and $t>4$. We first define the automaton
\texttt{triplegap} which recognizes these triples and then test if
the triples have the desired property.
\begin{center}
    \texttt{def triplegap "?msd\_fib En (t>0) \& (u>0) \& (v>0) \& \$marker2(n) \& \$marker2(n+t) \& \$marker2(n+t+u) \& \$marker2(n+t+u+v) \&\newline
(Aw (w>0) \& (w<t+u+v) \& (w!=t\&w!=t+u) => \~{}\$marker2(n+w))":}\smallbreak

\texttt{eval test "?msd\_fib At,u,v \$triplegap(t,u,v)\ <=>\newline ( \$isfib(t) \& \$isfib(u) \& \$isfib(v) \& \$isfib(t+v) \&\newline (t<=u) \& (u<=v) \& (t>4) )":}
\end{center}
\texttt{Walnut} says \texttt{TRUE}. The words $W_i$ between the markers are indeed
of the observed lengths. 
We now verify that these words are palindromes. 
Suppose that $W$ is marked, i.e., $\textcolor{red}2W\textcolor{red}2$, 
and that the markers correspond to the numbers $m<n$.
Then the digits in
$W$ correspond to numbers $k$ in
\href{https://www.oeis.org/a260317}{A260317} for $m+2<k<n$.
The word is a palindrome if and only if the reversal $x\mapsto n+m+2-x$ 
 of the
interval $[m+2,n]$
preserves this set of numbers $k$.

\begin{center}
    \texttt{eval test "?msd\_fib Ak,m,n\ ( m>14 \& k>m+2 \& n>k )\newline 
    (\& \$marker2(m) \& \$marker2(n) \& (Ap((p<n) \& (p>m)) => \~{}\$marker2(p)))
=> ( \$a260317(k) <=> \$a260317(n+m+2-k) )":
}
\end{center}
\texttt{Walnut} says \texttt{TRUE}. To complete the analysis, we need to show that each
factor $\textcolor{red}2U\textcolor{red}2V\textcolor{red}2W\textcolor{red}2$ such that
$|U|=|V|$ satisfies $U=V$. 
\begin{center}
\texttt{
eval testUV "?msd\_fib An,s
(n>14\&s>0\&\$marker2(n)\&\$marker2(n+s)\newline\&\$marker2(n+2*s)\&(At(t>0\&t<2*s\&t!=s)=>\~{}\$marker2(n+t)))
\newline=>(Ak (k>0\&k<s)=>(\$a260317(n+k)<=>\$a260317(n+s+k)))":
}
\end{center}
\texttt{Walnut} says \texttt{TRUE}. We check that $V$ coincides with $W$ for the first
$\lceil|W|/2\rceil$ digits from $W_3$ on,
which is the palindrome following the marker $\textcolor{red}{2}$ at number 40 in A260317.
\begin{center}
    \texttt{eval testVW "?msd\_fib An,s,u
(n>14\&s>0\&u>s\&\$marker2(n)\&\$marker2(n+s)\newline\&\$marker2(n+s+u)\&(At(t>0\&t<s+u\&t!=s)=>\~{}\$marker2(n+t)))
\newline=>(Ak (k>0\&k<=(u+1)/2)=>(\$a260317(n+k)<=>\$a260317(n+s+k)))":}
\end{center}
\texttt{Walnut} says \texttt{TRUE}. A260371 indeed has the observed properties.

In the introduction we saw that each prefix of length $F_n-2$ of $\mathbf F$ is a palindrome. Therefore, if you 
put a digit $3$ at the front and at the back of such a prefix, you get a palindrome of length $F_n$.
We just need to check that these palindromes are our $W_{n-1}$. Of course 
this can be done with \texttt{Walnut},
but we do it the old fashioned way. The prefixes $P_n$ of $\textbf F$ of Fibonacci length $F_{n+1}$
 can be defined inductively by $P_1=5$, $P_2=53$ and
$P_{n+1}=P_nP_{n-1}$. Removing the final two digits of $P_n$ and padding it with $3$'s gives
a palindrome $Q_n$ of length $F_{n+1}$. We have $Q_{n+1}=3P_nT_{n-1}3$ where $T_{n-1}$ is the
truncation of $P_{n-1}$ by its final two digits. The palindromes $Q_{n+1}=3P_{n-1}P_{n-2}T_{n-1}$
and $Q_n=3P_{n-1}T_{n-2}3$ coincide on the initial $F_n+F_{n-1}-1$ letters, which is more
than half of the length of $Q_{n+1}$ for $n\geq 4$. Therefore, $Q_4$ determines all subsequent
palindromes. It is equal to $W_4$. We conclude that the limit of the palindromes $W_n$ is equal
to~$3\mathbf F$.

\section{Polynomials without zero coefficients}

Fine~\cite{fine1947} proved that if all the coefficients of the polynomial $(x+1)^n$ are odd, then
$n=2^m-1$ for some $m$. We can verify that with~\texttt{Walnut}. These coefficients
are the binomials $\binom{n}{k}$ for $0\leq k\leq n$. Fine's result depends
on a classical theorem of Lucas, according to which for integers $k,\ n,$ and prime $p$, the following holds:
\[
\binom{n}{k}=\prod_{j=0}^r \binom{n_i}{k_i}\mod(p)
\]
where $k_i,n_i$ are the digits of the base $p$ expansions of $k=k_rp^r+\ldots+k_1p+ k_0$ and 
${n=n_rp^r+\ldots+n_1p+n_0}$, and we use the convention
that $\binom ab=0$ if ${a<b}$. In particular, $\binom nk$ is odd if and only if 
there does not exist an $i$ such that $k_i=1$ and $n_i=0$. In other words,
all digits $(k_i,n_i)$ are in $\{(0,0),(0,1),(1,1)\}$. Translating
this into \texttt{Walnut} syntax gives
\begin{center}
    \texttt{reg bincoef msd\_2 msd\_2 "([0,0]|[0,1]|[1,1])*":}
\end{center}
The binary numeration system is given by \texttt{msd\_2}.
A power of 2 in binary is a 1 followed by a string of 0's:
\begin{center}
    \texttt{reg power2 msd\_2 "0*10*":}
\end{center}
Fine's theorem in first order logic is
\[
\forall n\ \left(\exists m\ n+1=2^m\right)\Leftrightarrow \left(k\leq n\implies \binom{n}{k}=1\ \mathrm{mod} \ (2)\right).
\]
We can now check that with \texttt{Walnut}, where we write $n$ as $j$ to adhere to the convention that variables are entered in alphabetical order:
\begin{center}
    \texttt{eval fine "Aj \$power2(j+1) <=> (Ak (k<=j => \$bincoef(k,j)))":}
\end{center}
which evaluates as \texttt{TRUE}.

In a recent paper Hajdu et al~\cite{hajdu2023} extended Fine's theorem to polynomials modulo three.

\begin{theorem}[Hajdu et al]\label{Hajdu}
    If none of the coefficients of $(x-1)^c(x+1)^d$ is divisible by 3, then $c+d+1$ is
    of the shape $3^j, 2\cdot 3^j, 3^i+3^j, 2\cdot 3^i+3^j$ for $i>j\geq 0$
\end{theorem}

The question we would like to begin to address here can be formulated
as follows: what can be said about the degree $\deg f$ of polynomials
$f$ over a finite field $\mathbb{F}_q$ of $q$ elements for which all irreducible
factors are linear without zero coefficients?

The theorems of Fine and Hajdu et al show that for $q=2$ and $q=3$
it holds that the number of non-zero digits in the binary,
resp.~ternary expansion of
$(\deg f)+1$ is bounded, by 1 and 2 respectively.

As a special case of a deep result by Adamczewski et al~\cite{adamczewski2023} on formal power series,
it is possible to find an automaton with output that will help us
answer such a question. This automaton \texttt{co}$q$, will consist
of $q^{q-1}$ states.
The purpose of this automaton is that,
upon input a $q$-tuple $(k, a_1, a_2, \cdots, a_{q-1})$, it
outputs the $k$-th coefficient $\textrm{co}_q(k, a_1, ..., a_{q-1})$ of the polynomial
$$f=(x-\alpha_1)^{a_1}(x-\alpha_2)^{a_2}\cdots 
(x-\alpha_{q-1})^{a_{q-1}}\in\mathbb{F}_q[x],$$
where the $\alpha_i$ are the distinct elements of $\mathbb{F}_q$. The existence of such an automaton follows from the observation that,
summing over all tuples $0\leq  k, a_1, ..., a_{q-1} < \infty$,
$$\sum \textrm{co}_q(k, a_1, ..., a_{q-1}) x^k y_1^{a_1}\cdots y_{q-1}^{a_{q-1}} = \prod_{i=1}^q\frac{1}{(1-y_i(x-\alpha_i))},$$
(using the familiar geometric series identity repeatedly), which is
multivariate power series that is {\it rational} over $\mathbb{F}_q(x, y_1, \cdots, y_{q-1})$. Automaticity then results from the multivariate generalization of Christol's theorem by Adamczewski and Bell \cite{adamczewskibell2012}, while \cite{adamczewski2023} contains an
efficient algorithm to compute the automaton.
Using an implementation of this algorithm in Magma \cite{magma1997},
we generated the automata for $q=3, 4$, and $5$. They become
large very quickly, since one has to specify the $q^q$ possible
transitions for each of the $q^{q-1}$ states: leading to a specification
for a text file in \texttt{Walnut} with almost 2 million lines for $p=5$.

Starting with the modest automaton
\texttt{co3} having 9 states and 243 transitions, one computes any
coefficient of $(x-1)^c(x+1)^d$ modulo three from three \texttt{lsd\_3} 
integers $k, c, d$ as input. 

To verify its correctness, the effect of
multiplication by $x-1$ and $x+1$ for each state can be checked. 
This gives eighteen verifications in total; we
list three of them (as predicates without the \texttt{eval} command, to save space):

\vspace{0.1cm}
{\small
\begin{verbatim}
"?lsd_3 Aa,c,d(co3[a+1][c][d]=@0 & co3[a][c][d]=@0) => co3[a+1][c+1][d]=@0":
"?lsd_3 Aa,c,d(co3[a+1][c][d]=@1 & co3[a][c][d]=@0) => co3[a+1][c+1][d]=@2":
"?lsd_3 Aa,c,d(co3[a+1][c][d]=@2 & co3[a][c][d]=@0) => co3[a+1][c+1][d]=@1":                                              
\end{verbatim}
}

\normalsize\noindent
These check, for each of the three possible values,
that coefficient $a+1$ of 
$(x-1)\cdot f$ is obtained as the difference of coefficients
$a$ and $a+1$ of the original $f$.

\begin{figure}[t]
    \centering
    \includegraphics[width=0.7\linewidth]{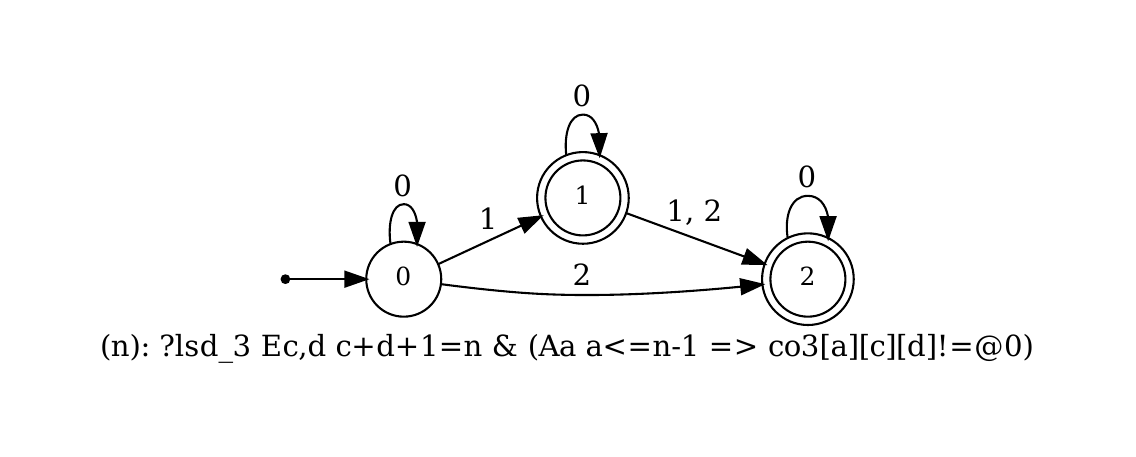}
    \caption{The {\texttt{lsd\_3}}-automaton that accepts $n=c+d+1$ such that
    none of the coefficients in $(X-1)^c(X+1)^d$ are divisible by three. Accepting
    states are marked by double circles.} 
    \label{fig:coef3}
\end{figure}

With \texttt{co3} it is possible in \texttt{Walnut} to check that Theorem~\ref{Hajdu} holds, as follows.
\vspace{0.1cm}
\newline
\texttt{ 
\hfill eval co3no0 "?lsd\_3 Ec,d c+d+1=n \& (Aa a<=n-1 => co3[a][c][d]!=@0)":}
\vspace{0.1cm}
\newline
The resulting automaton \texttt{co3no0} for \texttt{n} is shown 
in Fig.~\ref{fig:coef3}.
This does indeed confirm Theorem~\ref{Hajdu}, as the largest number
of non-zero digits in an input path (for $n$) leading to an accepting
state is equal to 3.

Hajdu et al ask if a far reaching generalization of these 
results holds: is it true that for every prime $p$ there exists
a constant $c_p$ such that any monic polynomial
$f(x)\in\mathbb Z[x]$ with only rational roots and no coefficients divisible by $p$ has at most $c_p$ nonzero digits in the base $p$
expansion of $\mathrm{deg}(f)+1$? In particular, is this true for $c_p=p-1$?
Much more will be said about this elsewhere, but using \texttt{Walnut}
we can prove two results in addition to the cases for $p=2$ and $p=3$ above.

The automaton \texttt{co5} we produced along the lines sketched 
above computes the coefficients of
$(X-2)^b(X-1)^c(X+1)^d(X+2)^e$ modulo five, but our \texttt{Walnut} computation for the
degrees of polynomials without coefficients divisible by five 
spirals out of control even on a powerful machine. 
As before, for \texttt{co3}, the correctness of \texttt{co5}
could inductively be checked easily (using $5^3$ checks in all).
And although it was also possible to verify the conjecture
for special cases (essentially leaving out some factors in
the product),
it was felt that the general case for $p=5$ ought to be doable.

\begin{figure}[t]
    \centering
    \includegraphics[width=0.7\linewidth]{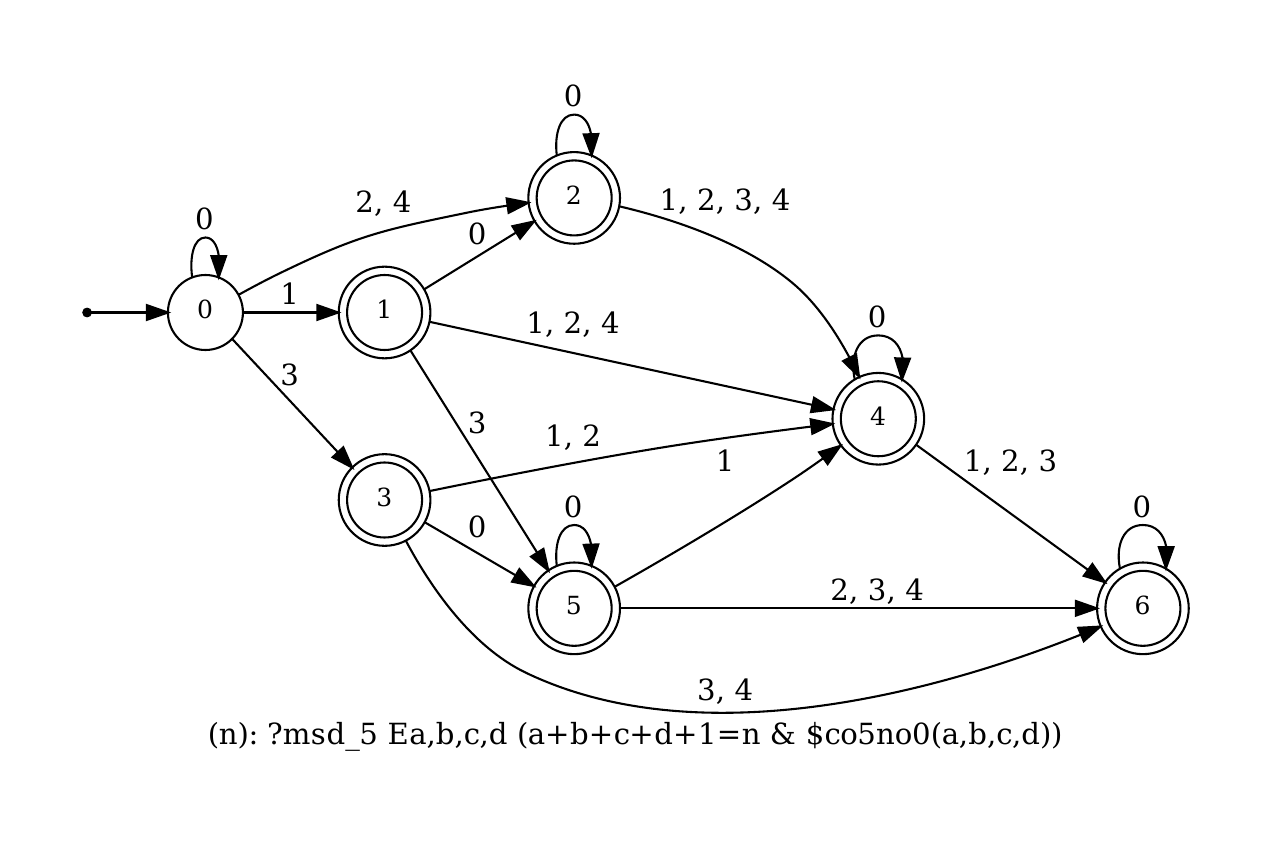}
    \caption{The {\texttt{msd\_5}}-automaton that accepts $n=b+c+d+e+1$ such that
    none of the coefficients in $(X-2)^b(X-1)^c(X+1)^d(X+2)^e$ is divisible by five. Accepting
    states are marked by double circles.} 
    \label{fig:lem5}
\end{figure}

Here John Nicol and Markus Frohme came to the rescue: they were
working on improvements for \texttt{Walnut} (version 7), to be reported
on in a forthcoming paper \cite{walnut7}. This version did succeed
in producing the analogon \texttt{co5no0} of \texttt{co3no0}.
The result is shown in Fig.~\ref{fig:lem5}.
Again, confirmation of the claims lies in the fact that no path to an
accepting state contains more than 4 non-zero symbols.

\begin{figure}[h]
    \centering
    \includegraphics[width=0.7\linewidth]{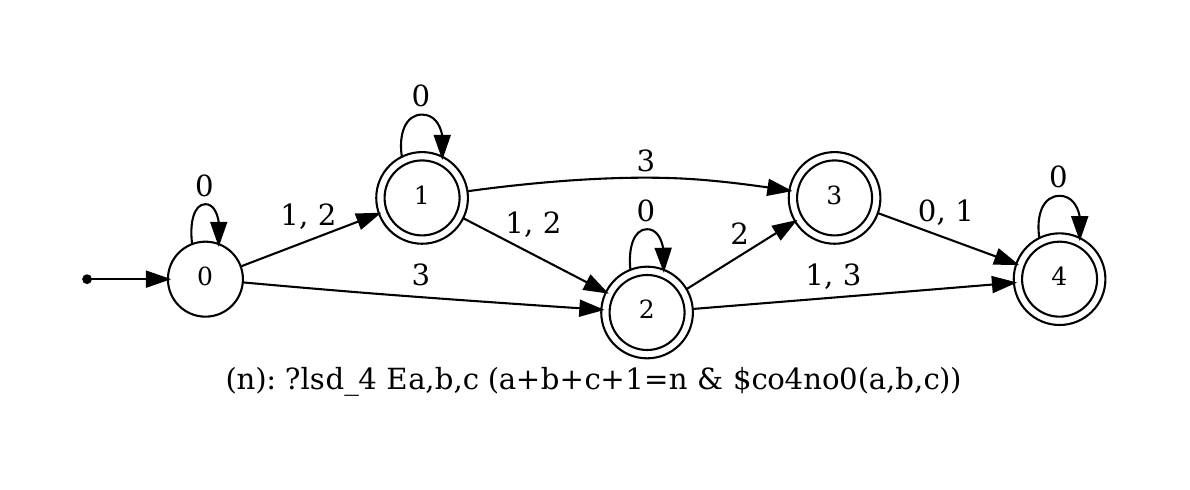}
    \caption{The {\texttt{lsd\_4}}-automaton that accepts $n=b+c+d+1$ such that
    none of the coefficients of product of linear factors of degree $n-1$ is zero. Accepting
    states are marked by double circles.} 
    \label{fig:lem4}
\end{figure}

The conjecture of Hajdu refers to arithmetic modulo primes $p$,
so in a prime field only. We decided to also consider the case
of $\mathbb{F}_4$, the field of 4 elements. With some care our
algorithm also produced the automaton \texttt{co4}, allowing
four inputs to produce the $k$-th coefficient of
$f=(x+1)^{b}(x+g)^{c}
(x+g^2)^{d}\in\mathbb{F}_4[x],$ where $g$ satisfies $g^2=g+1$ over
$\mathbb{F}_2$. After checking correctness, the automaton
\texttt{co4} easily produced \texttt{co4no0}, shown in 
Fig.~\ref{fig:lem4}.
Here we see that 4 is the largest number of non-zero digits in the
4-ary expansion of $(\deg f)+1$. As an example, the path $1, 1, 2, 1$
in the automaton, corresponding to $n=1+4+2\cdot 4^2+4^3=101$
suggests that there exist polynomials over $\mathbb{F}_4$ of degree
100 with no non-zero coefficients that are a product of linear factors.
Indeed,
$$f=(x+1)^{28}(x+g)^{35}(x+g^2)^{37}$$
is readily exhibited by Magma as an example. An automaton recognizing these
polynomials is implicitly also constructed by Walnut.


\section{Anti recurrence sequences}

Don't be unprepared for double negations!
Fibonaccis are generated by sums of numbers in the sequence. In contrast, \emph{anti-Fibonaccis} 
are sums of numbers that are \emph{not} in the sequence. Start with 0. The first two missing numbers are
1 and 2. Add them to get the next anti-Fibonacci~$3$. Continue with the next two missing numbers 4 and 5 to get
the anti-Fibonacci 9, and so on. This gives us sequence \href{https://oeis.org/A075326}{A075326}
\[
0,\ 3,\ 9,\ 13,\ 18,\ 23,\ 29,\ 33,\ 39,\ 43,\ 49,\ 53,\ \ldots 
\] 
The missing numbers are known as the \emph{non-anti-Fibonaccis} and
they can be found in \href{https://oeis.org/A249031}{A249031}:
\[
1,\ 2,\ 4,\ 5,\ 6,\ 7,\ 8,\ 10,\ 11,\ 12,\ 14,\ 15,\ 16,\ \ldots
\]
Thomas Zaslavsky completely cleared up the structure of the anti-Fibonaccis
and non-anti-Fibonaccis in an {unpublished note}~\cite{zaslavsky},
specifying observations of Doug Hofstadter in an earlier unpublished note~\cite{hofstadter}.
Zaslavsky proved that the gaps between anti-Fibonaccis come in pairs $\{6,4\}$ and $\{5,5\}$. If these
pairs are labelled $a$ and $b$ respectively, then they form the fixed point of the period-doubling 
morphism $a\mapsto ab,\ b\mapsto aa$. This sequence is implemented as \texttt{PD} in \texttt{Walnut}
with $a$ represented by zero and $b$ represented by one.
Zaslavsky proved that 
the sequence $X_n$ of anti-Fibonaccis satisfies $X_{2k+1}=3+10k$ and
$X_{2k+2}=8+10k+\mathtt{PD}[k]$.
The non-anti-Fibonaccis $Y_n$ can be described in a similar way,
if we write the indices modulo four:
\[
\begin{array}{lll}
Y_{4k}&=&5k, \\
Y_{4k+1}&=&5k+1, \\
Y_{4k+2}&=&5k+2, \\
Y_{4k+3}&=&5k+3+\mathtt{PD}[k].
\end{array}
\]
We can implement $X_n$ and $Y_n$
as synchronized sequences:
\begin{center}
    \texttt{def antifib "Ek n=2*k+1\&x=3+10*k | Ek n=2*k+2\&PD[k]=@1\&x=9+10*k\newline | Ek n=2*k+2\&PD[k]=@0\&x=8+10*k":}

    \texttt{def nonafib "Ek n=4*k+1\&x=5*k+1 | Ek n=4*k+3\&PD[k]=@1\&x=5*k+4 |\newline Ek n=4*k+3\&PD[k]=@0\&x=5*k+3 | Ek n=4*k\&x=5*k | Ek n=4*k+2\&x=5*k+2":}
\end{center}
We now verify Zaslavsky's results that these two sequences are indeed the anti-Fibonaccis
and the non-anti-Fibonaccis (or missing numbers). 
The defining properties of these sequences are complementary and that each anti-Fibonacci is a sum of two
consecutive missing numbers $X_n=Y_{2n-1}+Y_{2n}$.
We first verify this additive relation between the two sequences:
\begin{center}
    \texttt{eval zaslavsky1 "Ak,x,y,z ((k>0)\&\$nonafib(2*k-1,x)\&\$nonafib(2*k,y)\&\$antifib(k,z))=>x+y=z":}
\end{center}
\texttt{Walnut} says \texttt{TRUE}. We now verify
that the sequences are complementary in two steps. Disjointness is checked by
\begin{center}
    \texttt{eval zaslavsky2 "Ei,j,n (n>0)\&\$nonafib(i,n)\&\$antifib(j,n)":}
\end{center}
on which \texttt{Walnut} returns \texttt{FALSE}. We verify that each
number is in one of the two sequences
\begin{center}
    \texttt{eval zaslavsky3 "An (n>0) => (Ej \$nonafib(j,n)|\$antifib(j,n))":}
\end{center}
on which \texttt{Walnut} returns \texttt{TRUE}. The sequences are indeed complementary,
establishing that the $X_n$ are the anti-Fibonaccis and the $Y_n$ are the non-anti-Fibonaccis.
\medbreak
Kimberling and Moser~\cite{kimberling2019linear} observed that 
there is another way to describe the non-anti-Fibonaccis $Y_n$
using the \texttt{mex} or minimal excluded value. For any set
of natural numbers the \texttt{mex} is the minimal element of its complement.
Create three sequences $A_n,B_n,C_n$ starting from $A_1=1, B_1=2, C_1=3$
such that $A_{n+1}=\mathtt{mex}\{A_i,B_i,C_i\colon i\leq n\}$ and
$B_{n+1}=\mathtt{mex}\left(\{A_i,B_i,C_i\colon i\leq n\}\cup\{A_{n+1}\}\right)$
and $C_{n+1}=A_{n+1}+B_{n+1}$. The anti-Fibonaccis are $C_n$. The odd-indexed
non-anti-Fibonaccis are $A_n$ and the even-indexed are $B_n$.
Instead of creating three sequences we can create more.

\begin{figure}[h]
    \centering
    \includegraphics[width=0.8\linewidth]{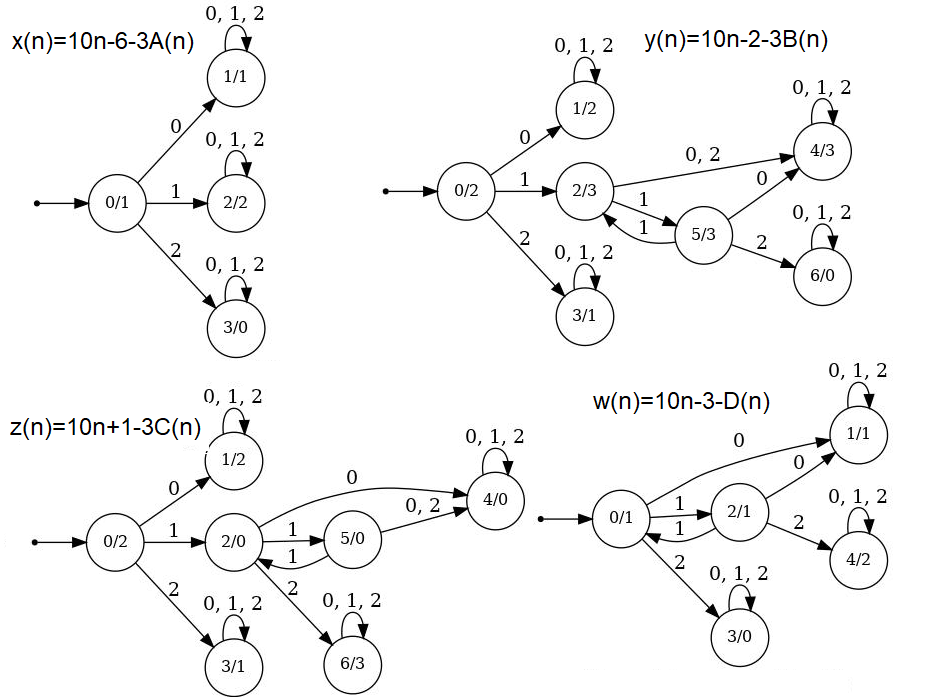}
    \caption{The guessed {\texttt{lsd\_3}}-automata} for the remainders associated to the anti-Tribonacci
     sequence.
    \label{fig:automata}
\end{figure}

If we apply the same process to four sequences $A_n,B_n,C_n,D_n$ 
instead of three, then we get the \emph{anti-Tribonacci} sequence.
As in the previous case, $A_n, B_n, C_n$ are defined by the \texttt{mex} and $D_n$ 
is the sum. The $D_n$ are the anti-Tribonaccis \href{https://oeis.org/A265389}{A265389}:
\[6,\ 16,\ 27,\ 36,\ 46,\ 57,\ 66,\ 75,\ 87,\ 96,\ \ldots\]
The other three sequences are \href{https://oeis.org/A297464}{A297464-6}.
Clark Kimberling conjectured that
\[
\begin{array}{rccl}
0 & \leq & 10n - 6 - 3A_n& \leq 2, \\
0 & \leq & 10n - 2 - 3B_n& \leq 3, \\
0 & \leq & 10n + 1 - 3C_n& \leq 3, \\
0 & \leq & 10n - 3\ - D_n &\leq 2.
\end{array}
\]
Following Zaslavsky's lead, we guess that these four equations form automatic
sequences and verify this using \texttt{Walnut}.
We call them the \emph{remainders} and our guessed
3-automata are given in Figure~\ref{fig:automata}.

\medskip

We have implemented
these as \texttt{lsd\_3} automata \texttt{xkimber, ykimber, zkimber, wkimber}.
The sequences $A,B,C,D$ are given by
\medskip

\texttt{
\hspace*{2.3em}def seqa "?lsd\_3 s=(10*n-6-xkimber[n-1])/3": \\
\hspace*{4em}def seqb "?lsd\_3 s=(10*n-2-ykimber[n-1])/3": \\
\hspace*{4em}def seqc "?lsd\_3 s=(10*n+1-zkimber[n-1])/3": \\
\hspace*{4em}def seqd "?lsd\_3 s=(10*n-3-wkimber[n-1])":
}

\noindent 
To prove Kimberling's conjecture we need to verify that these sequences
are complementary, defined by the \texttt{mex}, and that $A+B+C=D$. We first check that
$A_n<B_n<C_n<A_{n+1}$, proving that $A,B,C$ are disjoint:
\medskip

\texttt{
\hspace*{1em}eval test1 "?lsd\_3 An,s,t,u,v (\$seqa(n,s) \& \$seqb(n,t)\newline
\hspace*{2.2em} \& \$seqc(n,u) \& \$seqa(n+1,v)) => ((s<t)\&(t<u)\&(u<v))":
}
\medskip

\noindent
Next, we verify that the union of the four sequences is equal to $\mathbb N$:
\medskip

\texttt{
eval test2 "?lsd\_3 As (s>0)=>Em
(\$seqa(m,s)|\$seqb(m,s)|\newline 
\hspace*{18.2em}\$seqc(m,s)|\$seqd(m,s))":}

\medbreak
\noindent
In both cases \texttt{Walnut} says \texttt{TRUE}. We now verify that $D$
is disjoint from $A\cup B\cup C$ to establish that the sequences are complementary:
\medskip

\texttt{
eval test3 "?lsd\_3 Em,n,s \$seqd(n,s) \& \newline \hspace*{14.2em}(\$seqa(m,s)|\$seqb(m,s)|\$seqc(m,s))":
}
\medskip

\noindent \texttt{Walnut} returns \texttt{FALSE} which confirms that the sequences
are complementary. Finally, we check that $A+B+C=D$
\medbreak

\texttt{
eval test4 "?lsd\_3 An,s,t,u,v (\$seqa(n,s) \&
\$seqb(n,t) \&\newline\hspace*{11em} \$seqc(n,u) \& \$seqd(n,v)) => v=s+t+u":
}
\medbreak

\noindent \texttt{Walnut} says \texttt{TRUE}. 
From $A+B+C=D$ we easily deduce that $D_n>C_{n+1}$ and therefore $A,B,C$ are indeed defined
by the \texttt{mex}, confirming the conjecture. Kimberling also conjectured that
\[
A_{n+4}+A_{n}=A_{n+3}+A_{n+1}
\]
This is now a straightforward verification
\medbreak
\texttt{\hspace*{3em}eval kimconj "?lsd\_3 An,r,s,t,u (\$seqa(n,r)\\
\hspace*{4em}\&\$seqa(n+1,s)\&\$seqa(n+3,t)\&\$seqa(n+4,u))=>r+u=s+t":}

\medbreak

As any clergyman will tell you, there is no need to stop at four.
Indeed, Kimberling defined five complementary sequences
 $A_n,B_n,C_n,D_n,E_n$ in which the
$E_n$ are the \emph{anti-Teranaccis} and the other sequences are
the non-anti-Teranaccis (or missing numbers) defined by the \texttt{mex}.
Again, these sequences are complementary and the anti-Teranaccis are
sums of the missing numbers.
The anti-Teranaccis form sequence \href{https://oeis.org/A299409}{A299409}
where it is conjectured that
\[
\begin{array}{rccl}
0 & \leq & 17n - 11 - 4A_n& \leq 4, \\
0 & \leq & 17n -\ 7 - 4B_n& \leq 4, \\
0 & \leq & 17n -\ 3 - 4C_n& \leq \textcolor{red}3, \\
0 & \leq & 17n +\ 1 - 4D_n& \leq 3, \\
0 & \leq & 17n -\ 5 -\ E_n &\leq 3.
\end{array}
\]
Again, we can guess automata for the remainder sequences, which we call
\texttt{xx}, \texttt{yy}, \texttt{zz}, \texttt{vv}, \texttt{ww} this time, and define the sequences \texttt{A,B,C,D,E}
from them. It turns out that there is a typo in the conjecture.
The remainder of $C$ is bounded by 4 instead of \textcolor{red}{3}.
\begin{figure}
    \centering
    \includegraphics[width=0.5\linewidth]{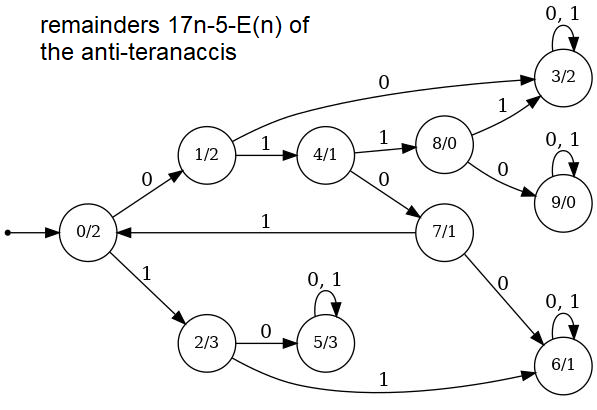}
    \caption{The guessed {\texttt{lsd\_2}-automaton}} for
    the anti-Teranacci remainder.
    \label{fig:2}
\end{figure}

\medbreak
\texttt{
\hspace*{2.3em}def seq4A "?lsd\_2 s=(17*n-11-xxkimber[n-1])/4": \\
\hspace*{4em}def seq4B "?lsd\_2 s=(17*n-7-yykimber[n-1])/4": \\
\hspace*{4em}def seq4C "?lsd\_2 s=(17*n-3-zzkimber[n-1])/4": \\
\hspace*{4em}def seq4D "?lsd\_2 s=(17*n+1-vvkimber[n-1])/4": \\
\hspace*{4em}def seq4E "?lsd\_2 s=(17*n-5-wwkimber[n-1])":
}
\smallbreak
\begin{figure}[h]
    \centering
    \includegraphics[width=1\linewidth]{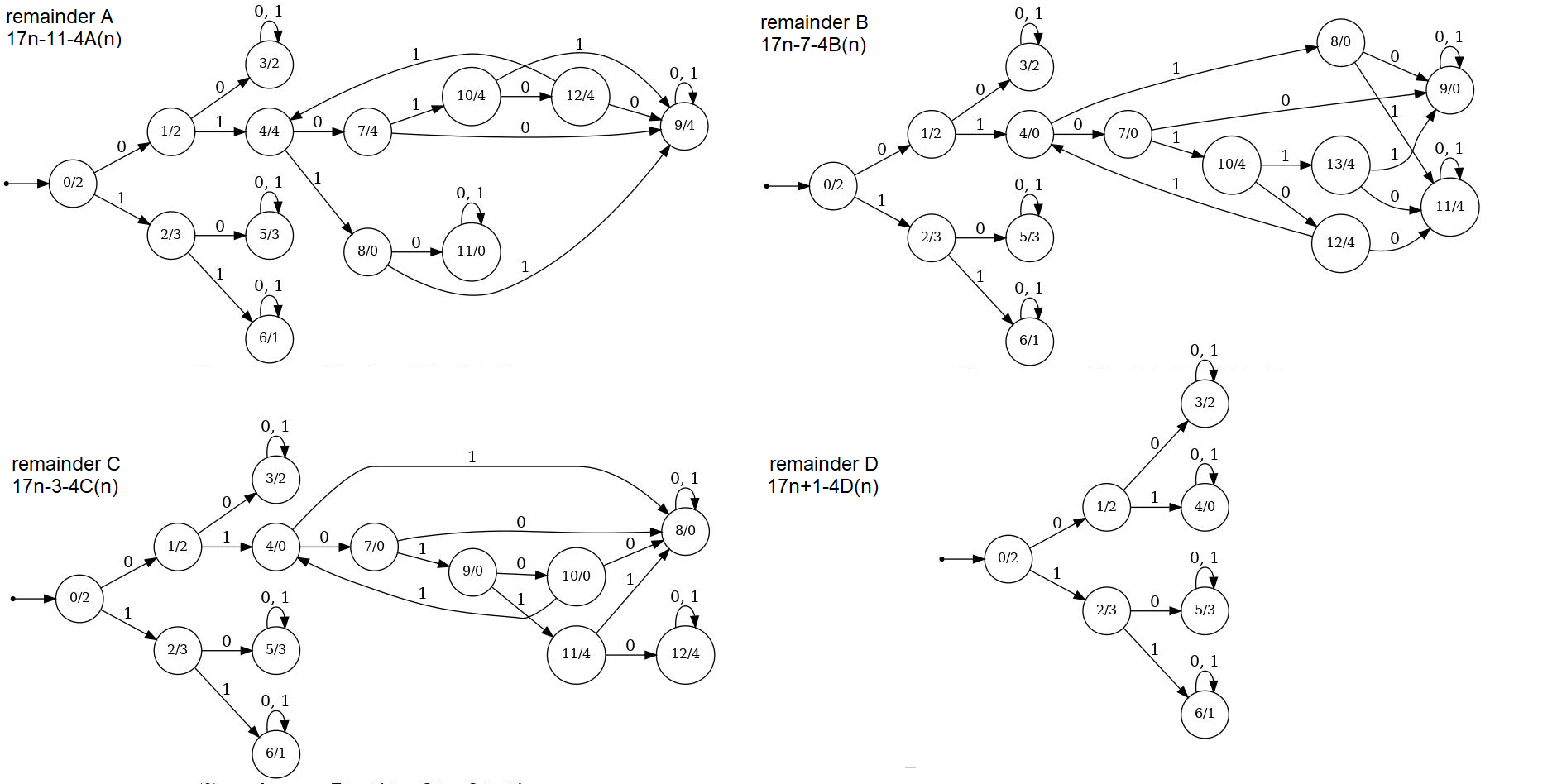}
    \caption{The four guessed {\texttt{lsd\_2}-automata}} for
    the remainders of the non-anti-Teranaccis.
    \label{fig:3}
\end{figure}

\medbreak
We run the same four tests which again return \texttt{TRUE, TRUE, FALSE, TRUE}
as required.
\medbreak
\texttt{\hspace*{0.5em}eval test1 "?lsd\_2 An,s,t,u,v,w (\$seq4A(n,s) \& \$seq4B(n,t)\\
\hspace*{4em}\& \$seq4C(n,u) \& \$seq4D(n,v) \& \$seq4A(n+1,w)) \\  \hspace*{4em}=> ((s<t)\&(t<u)\&(u<v)\&(v<w))":}
\smallbreak
\texttt{
eval test2 "?lsd\_2 As (s>0)=>Em (\$seq4A(m,s)|\$seq4B(m,s)|\\
\hspace*{4em}\$seq4C(m,s)|\$seq4D(m,s)|\$seq4E(m,s))": 
}

\medbreak
\texttt{
eval test3 "?lsd\_2 Em,n,s \$seq4E(n,s) \&\\
\hspace*{3.5em}
(\$seq4A(m,s)|\$seq4B(m,s)|\$seq4C(m,s)|\$seq4D(m,s))":
}
\smallbreak
\texttt{
eval test4 "?lsd\_2 An,s,t,u,v,w (\$seq4A(n,s) \& \$seq4B(n,t) \&\\
\hspace*{3.5em}\$seq4C(n,u) \& \$seq4D(n,v) \& \$seq4E(n,w)) => w=s+t+u+v":
}
\medbreak
We define the anti-$k$-naccis as the sums of
$k$ consecutive missing numbers. There is an obvious pattern.
\begin{conjecture}[The Clergyman's Conjecture.] Let $X_n$ be
the sequence of anti-$k$-naccis. The difference
$
X_n-(k^2+1)n
$
is $k$-automatic.
\end{conjecture}

\section{Subsumfreeness}
In a very similar fashion we consider general greedy 3-subsumfree 
sequences in this section. 
Consider a sequence $a$ with start values $a_1, a_2, a_3$ and further
values $a_n$, for $n>3$ satisfying that $a_n$ is the smallest number
exceeding $a_{n-1}$ which is not the sum of three different
previous entries: so $a_n\neq a_i + a_j + a_k$ for
$1\le i < j < k < n$. In his note \cite{stephan2004} Stephan lists 
several conjectures (in Section 2.3.1) of the following kind.
If $a=(1,2,3,\ldots)$ then
$a_{n+6}-a_{n+5}=a_{n+1}-a_n$ for $n>6$.

We will use the notation $S_{x,y,z}$ for the greedy 3-subsumfree sequence with starting values $x<y<z$. The cases $S_{1,2,3}$ and $S_{1,3,4}$ appear as
\href{https://www.oeis.org/a026471}{A026471} and \href{https://www.oeis.org/a026475}{A026475} in the OEIS.

Since it was not clear to us how to verify such a conjecture using \texttt{Walnut}, we first reworked them in a more explicit form that is amenable to direct verification, and then attempted to do this in \texttt{Walnut}. This led to families of conjectures of increasing generality (see Conjectures \ref{conj123}, \ref{conj_g} and \ref{conj_d_g} below) for (almost all) greedy 3-subsumfree sequences starting with 1. The example for $(1,2,3,\ldots)$ reads as follows. 
\begin{conjecture}\label{conj123}
Let $z$ be a positive integer. Then
$$z\in S_{1,2,3}\quad\iff\quad z\in{1,5,13}\textrm{\ \ or\ \ }z\bmod 23\in\{2,3,4,14,15\}.$$
\end{conjecture}
It is easy to verify by hand that $S_{1,2,3}$ starts as:
$$S_{1,2,3}=1,2,3,4,5,13,14,15,25,26,27,37,38,48,49,50,60, 61, \ldots$$
and the conjecture would imply that from the seventh entry `$14$' on the sequence modulo $23$ is periodic with period 5.

With a few lines of \texttt{Walnut} code we can indeed verify this! First define five relevant residue classes modulo 23: 

\medskip
\texttt{
\hspace*{2.3em}def is02 "Ek z=k*23+2": \\
\hspace*{4em}def is03 "Ek z=k*23+3": \\
\hspace*{4em}def is04 "Ek z=k*23+4": \\
\hspace*{4em}def is14 "Ek z=k*23+14": \\
\hspace*{4em}def is15 "Ek z=k*23+15":
}
\medbreak\noindent
which enables us to define the characteristic function for the sequence from the conjecture:

\medbreak
\texttt{
\hspace*{0.5em}def seq123 "(z=1 | z=5 | z=13 | \\
\hspace*{4em} \$is02(z) | \$is03(z) | \$is04(z) | \$is14(z) | \$is15(z) )":
}
\medbreak\noindent
and we can now test the property that a positive integer $z$ is in the sequence if and only if if it is not the sum of 3 previous entries:
\medbreak
\texttt{
\hspace*{0.5em}eval prop "Az z>0 => \\
\hspace*{3em}(\$seq123(z) <=> z>0 \& $\sim$(E a, b, c a<b \& b<c \& \\
\hspace*{8em}\$seq123(a) \& \$seq123(b) \& \$seq123(c) \& a+b+c=z))"::}
\medbreak\noindent
\texttt{Walnut} returns \texttt{TRUE} in little over a second: our conjecture is now a Theorem! We state (and proved by ) this also for two more cases, $S_{1,3,4}$ and $S_{1,4,5}$.
\begin{theorem}For positive integers $z$:
\begin{eqnarray*}z\in S_{1,2,3}&\iff& z\in{1,5,13}\textrm{\ \ or\ \ }z\bmod 23\in\{2,3,4,14,15\};\\
z\in S_{1,3,4}&\iff& z\in{1,7,19}\textrm{\ \ or\ \ }z\bmod 33\in\{3,4,5,6,20,21,22\};\\
z\in S_{1,4,5}&\iff& z\in{1,9,25}\textrm{\ \ or\ \ }z\bmod 43\in\{4,5,6,7,8,26,27,28,29\}.
\end{eqnarray*}
\end{theorem}
\noindent
%

In fact we believe the following generalized conjecture to hold.

\begin{conjecture}\label{conj_g}
For every $g\geq2$ the greedy 3-subsumfree sequence $S_{1,g,g+1}$ is
characterized as follows: 
\begin{eqnarray*}
z\in S_{1,g,g+1}&\iff& z\in\{1,2g+1,6g+1\}\textrm{\ \ or\ \ }z\bmod 10g+3\in\\
&&\{g, g+1,\cdots, 2g\}\cup\{6g+2, 6g+3,\cdots,7g+1\}.\end{eqnarray*}
In particular, after the first $g+4$ entries the sequence modulo $10g+3$ is periodic with period $2g+1$.
\end{conjecture}
Such parametrized statement is not so easily proved in \texttt{Walnut}. To be more precise, what one could attempt is to define the following:
\medbreak
\texttt{
\hspace*{0.5em}def isinGmodm "Ek,w ( ((w>=g \& w<=2*g)|(w>=6*g+2 \& w<=7*g+1))\\
\hspace*{14em} \& z=k*(10*g+3)+w)":\\
\hspace*{2.1em}def isinG "z=1 | z=2*g+1 | z=6*g+1 | \$isinGmodm(g,z)":}
\medbreak\noindent
and then universally quantify (over $g$) this property:
\medbreak
\texttt{
\hspace*{0.5em}eval prop "Az z>=g+4 => \\
\hspace*{3em}(\$isinG(z) <=> $\sim$(E a, b, c a<b \& b<c \& a+b+c=z \& \\
\hspace*{12em} \$isinG(a) \& \$isinG(b) \& \$isinG(c) ))":}
\medbreak\noindent
The problem here is that the definition of \texttt{isinGmodm} is not allowed in \texttt{Walnut}
\medbreak
\texttt{
\hspace*{0.5em}the operator * cannot be applied to two variables}
\medbreak\noindent
Whenever a numerical value is substituted for $g$ in the definition
of \texttt{isinGmodm} and \texttt{isinG} all is fine. In fact we did this for all $g\in\{2,3,\cdots 10\}$ successfully.

As a matter of fact we have firm 
computational evidence, from an implementation in Magma \cite{magma1997}, for the following meta-Conjecture.
\begin{conjecture}\label{conj_d_g}
Let $d\geq 2$. For every $g\geq d+1$ the greedy 3-subsumfree sequence $S_{1,g,g+d}$ is
characterized as follows:
\begin{eqnarray*}
z\in S_{1,g,g+d}&\iff& z\in\{1,g,2g+d-1,2g+d\}\textrm{\ \ or\ \ }z>g+d \textrm{\ and\ }\\
&& z\bmod 5g+2d\in\{g+d-2, g+d-1,\cdots, 2g+d-2\}.\end{eqnarray*}
In particular, for $d\geq 2$ and every $g\geq d+1$ after the first $g+3$ entries (in a preperiod) the sequence $S_{1,g,g+d}$ modulo $5g+2d$ is periodic with period $g+1$.
\end{conjecture}
\begin{remarks}\rm
Note that Conjecture \ref{conj_g}, although of the same form,
is not a special case of the 
meta-conjecture, as the specific values (for modulus and
(pre)period length) are irregular.

Also note that most, but not {\it all}, greedy 3-subsumfree sequences are
covered by the meta-conjecture: usually the cases of small values for
$g$ for given $d$ are special in the sense that the indicated
modular periodicity occurs for deviating values of modulus $m$
and period length $p$. Below is one example.

In the regular cases of the meta-conjecture, it is also possible
to state explicitly what the shape of the period (and preperiod) will be.
\end{remarks}
\begin{example}\rm
For $d=7$ and $g=5$ the sequence $S_{1,5,12}$ becomes
periodic modulo 321 with period length 32.
\end{example}
Here, by way of example, is a verification in  of the case $d=4, g=4$ of
Conjecture \ref{conj_d_g}. We first define the right residue classes modulo $5g+2d=28$
and then ask whether or not the 3-subsumfree property holds for all entries in $S_{1,4,8}$
after $1, 4, 8$:
\medbreak
\texttt{
\hspace*{0.5em}def isres "Eh, k h>= 6 \& h<=10 \& z=k*28+h":\\
\hspace*{2.1em}def seq148 "(z=1|z=4|z=11|z=12|(z>=8 \& \$isres(z)))":\\
\hspace*{2.1em}eval prop "Az z>8 => (\$seq148(z) <=> $\sim$(E a, b, c a<b \& b<c \\
\hspace*{4.1em}\& \$seq148(a) \& \$seq148(b) \& \$seq148(c) \& a+b+c=z))"::}
\medbreak\noindent
which returns {\tt TRUE} almost instantly.

Finally, it may be worthwhile to say something 
about the manner in
which these conjectures were obtained. At first sight it is
not clear at all that the definition of subsumfree sequence
leads to automatic sequences. Since this is important for the
success of any attempt to invoke \texttt{Walnut}, it is useful to see how
it arises naturally. But also, although not fitting in with the 
nature of this paper, it is useful for finding `pen and paper' proofs
for the claims, some of which we hope to present elsewhere.

The simple case is that it pays off to look at the sequence of
first differences of these subsumfree sequences. Consider again
$S_{1,2,3}$. A straightforward computation leads to this initial segment
and the corresponding sequence of differences:
$$S_{1,2,3}=1, 2, 3, 4, 5, 13, 14, 15, 25, 26, 27, 37, 38, 48, 49, 50, 60, 61, 71, 72, 73, \cdots$$
and
$$D=1, 1, 1, 1, 8, 1, 1, 10, 1, 1, 10, 1, 10, 1, 1, 10, 1, 10, 1, 1, \cdots .$$
Considering this it is not difficult to
infer the conjectured periodicity.
\section{Conclusion}
 
We would like to thank John Nicol and Markus Frohme for running the
\texttt{co5} automaton on their upcoming new version of \texttt{Walnut}, and 
Jeffrey Shallit for helping us out with \texttt{Walnut} commands and
for giving a lecture in our course.

\bibliographystyle{siam}
\bibliography{deBE}

\end{document}